\documentclass[USenglish]{article}	

\usepackage[utf8]{inputenc}				
\usepackage[big,online]{dgruyter}	
\usepackage{lmodern}
\usepackage{microtype}
\usepackage[numbers,square,sort&compress]{natbib}

\theoremstyle{dgthm}
\newtheorem{theorem}{Theorem}

\newtheorem{proposition}{Proposition}

\theoremstyle{dgdef}
\newtheorem{definition}{Definition}

\newtheorem{remark}{Remark}

\usepackage{pstricks}
\usepackage{amsfonts}
\usepackage{dsfont}
\newcommand{\N}{\mathbb{N}}

\newcommand{\R}{\mathbb{R}}
\usepackage{etoolbox}%
\AtEndEnvironment{proof}{\null\hfill\qedsymbol}%
\begin{document}

\journalname{Random~Oper.~Stoch.~Equ.}
\DOI{10.1515/rose-XXXX}
\articletype{Research Article}
\journalyear{2022}
\journalvolume{30}
\journalissue{1}
\startpage{1}

\communicated{A. F. Turbin}

\title{Reflected BSDEs with two completely separated barriers and regulated trajectories in general filtration.}

\runningtitle{RBSDEs with two completely separated barriers and regulated trajectories.}

\author*[1]{Brahim BAADI}

\runningauthor{B. BAADI.}
\affil[1]{\protect\raggedright
Ibn Tofa\"{\i}l University, Department of mathematics, Faculty of sciences, K\'{e}nitra, Morocco, e-mail: brahim.baadi@uit.ac.ma}



\received{April 7,  2021}
 \accepted{December 4, 2021}

\abstract{In this paper, we study doubly reflected  Backward Stochastic Differential Equations defined on probability spaces equipped with filtration satisfying only the usual assumptions of right continuity and completeness in the case where the barriers $L$ and $U$ are not necessarily  right continuous. We suppose that the barriers $L$ and $U$ and their left limits are completely separated and we show existence and uniqueness of the solution.}

\maketitle
	
\section{Introduction.}
In this paper, we study  the problem of existence and uniqueness of the solution of backward stochastic differential equations (BSDE) with two reflecting optional barriers (or obstacles) $L$ and $U$.  Our main in this work is to deal with equations on probability space with general filtration $\mathbb{F}=\{\mathcal{F}_{t}, t\geq 0\}$ satisfying only usual conditions of right continuity and completeness. Also, we assume  that the lower barrier L and the upper barrier U are completely separated in the sense that ($L_{t}<U_{t}$) and ($L_{t-}<U_{t-}$) for all $t\in[0,T]$  and which  are two regulated process, i.e. processes whose trajectories have left and right finite limit. Consequently, the solution of these equations need not  be c\`{a}dl\`{a}g but are called regulated processes.

Precisely, a solution for the BSDE with two reflecting barriers associated with a generator $f(t,y)$, a terminal value $\xi$, a lower barrier $L$ and an upper barrier $U$   (RBSDE$(\xi,f,L,U)$ for short), is a quadruple of processes $(Y,M,K,A)$ which mainly satisfies:
$$
 \left\{
   \begin{array}{ll}
    Y_{t}=\xi+\int_{t}^{T}f(s,Y_{s})ds+(K_{T}-K_{t})-(A_{T}-A_{t})-\int_{t}^{T}dM_{s}, t\in[0,T]. \\
     L_{t}\leq Y_{t}\leq U_{t}, t\in[0,T], and  \\
     \int_{0}^{T}(U_{s-}-Y_{s-})dA_{s}^{*}+\sum_{s<T}(U_{s}-Y_{s})\Delta^{+}A_{s}=0, and \\
     \int_{0}^{T}(Y_{s-}-L_{s-})dK_{s}^{*}+\sum_{s<T}(Y_{s}-L_{s})\Delta^{+}K_{s}=0.
\end{array}
 \right.
$$
where $Y$ has regulated trajectories, $K$, $A$ are increasing processes
such that $K_{0}=A_{0}=0$, $M$ is a local martingale with $M_{0}=0$, $K^{*}$ (resp. $A^{*}$) the c\`{a}dl\`{a}g part of $K$  (resp. $A$) and  $\Delta^{+}K$ (resp. $\Delta^{+}A$) the
right jump of $K$ (resp. of $A$). The reason we chose the minimality conditions
$$
\int_{0}^{T}(Y_{s-}-L_{s-})dK_{s}^{*}+\sum_{s<T}(Y_{s}-L_{s})\Delta^{+}K_{s}=0
\,\,\,\ \text{and} \,\,\,\
\int_{0}^{T}(U_{s-}-Y_{s-})dA_{s}^{*}+\sum_{s<T}(U_{s}-Y_{s})\Delta^{+}A_{s}=0,
$$
is to use the penalization method for regulated BSDE with regulated trajectories proposed  by Klimsiak et al. in \cite{Klimsiaketal2016}.  Note that if  $L$ and $U$ are c\`{a}dl\`{a}g, then  this condition reduces  to the classical condition: \cite[(1.3)]{Klimsiaketal2016}.

Generally speaking, in BSDE theory, during several years, there have been a lot of works which study the problem of existence and uniqueness of BSDE with two reflecting barriers under these  three conditions:
\begin{description}
  \item[a)]  one of the obstacles is  a semimartingale.
  \item[b)]  the  Mokobodski condition: between $U$ and $L$ one can find a process $X$ such that $X$ is a difference of nonnegative c\`{a}dl\`{a}g supermartinagles.
  \item[c)]  the barriers are completely separated: $L_{t}<U_{t}$ and $L_{t-}<U_{t-}$ for
all $t\in[0,T]$ a.s.
\end{description}
Under the assumption \textbf{b}), the problem is studied in \cite{Klimsiak2015},\cite{Klimsiak2013},\cite{CvitKara1996},\cite{BahlaHamadMez2005},\cite{CrepeyMatoussi2008}... in the case of continuous  or right-continuous  obstacles and/or a larger filtration than the Brownian, but the issue with this condition is that it  is quite difficult to check in practice. Then, it has been removed by Hamad\`{e}ne and Hassani in \cite{HamadeneHassani2005}, when they showed that if the assumption \textbf{c}) hold, the two barriers reflected BSDE has a unique solution. Under  the same assumption there are also a lot of works which dealt with the  problem of existence and uniqueness, for instance, the papers \cite{HamadeneHassaniOuknine2010},\cite{HamadeneHang2009},\cite{Hassairi2016},\cite{MateuszTopolewski2016}...
In all of the above-mentioned works (and others) on double reflected BSDEs, the barriers are
assumed to be at least right-continuous.

The only paper dealing with BSDEs with two reflecting  barriers that are not c\`{a}dl\`{a}g, in our knowledge, is the paper by Grigorova et al. in \cite{Grigorova2017}. The authors proved the existence and uniqueness of the solution of  double reflected BSDE with two irregular barriers  satisfying the generalized Mokobodzki's condition. First they showed the existence and uniqueness in the case where the driver does not depend on solution, then they proved a priori estimates for the doubly reflected BSDE by using Gal'chuk-Lenglart's formula and from these they derived the existence and uniqueness of the solution with general Lipschitz driver by using the  Banach fixed point theorem.

BSDE with two reflecting barriers have been introduced by Cvitanic and Karatzas
in \cite{CvitKara1996} in the case of continuous barriers and a Brownian filtration. The solutions of such equations are constrained to stay between two adapted barriers $L$ and $U$ with $L\leq U$ and $L_{T}= U_{T}$. In the case of the continuous/c\`{a}dl\`{a}g barriers, reflected doubly BSDE have been studied by several authors in \cite{HamadeneHassaniOuknine2010},\cite{HamadeneHang2009},\cite{Hassairi2016},\cite{MateuszTopolewski2016},
\cite{HamadeneHassani2006},\cite{HamadeneHdhiri2006}, \cite{Klimsiak2015},\cite{Klimsiak2013},\cite{CvitKara1996},\cite{BahlaHamadMez2005},
\cite{CrepeyMatoussi2008},\cite{DumQueSul2016},\cite{EssOukHarr2005},\cite{HamadeneLepeltier2000} and  \cite{Grigorova2017} (for  regulated  barriers case).

This paper is organized as follows:

In the second and third section, we give some preliminary and some result related to BSDE with one barrier (definition, existence). In section four, we recall the doubly reflected BSDE definition and we prove a comparison and uniqueness result. In the fifth section, we deal with the notion of local solution of doubly reflected BSDE, which is a solution  of that equation but between  two comparable stopping times. Some local solution properties are also  given. Section six is reserved to our main result of this paper.

\section{Preliminaries.}
Let us consider a filtered probability space $(\Omega,
\mathcal{F},\mathbb{P},\mathbb{F}=\{\mathcal{F}_{t}, t\geq 0\})$.
The filtration is assumed to be complete, right continuous and
quasi-left continuous.

Let $T>0$ be a fixed positive real number. We recall that a function $y:[0,T]\rightarrow \R^{d}$ is called regulated if for every $t\in[0,T]$ the limit $y_{t+}=\lim_{u\downarrow t}y_{u}$ exists, and for every $t\in[0,T]$ the limit $y_{t-}=\lim_{u\uparrow t}y_{u}$ exists. For any regulated function
$y$ on $[0,T]$, we denote by $\Delta^{+}y_{t}=y_{t+}-y_{t}$ the size of the
right jump of $y$ at $t$, and by $\Delta^{-}y_{t}=y_{t}-y_{t-}$ the size of the left jump of
$y$ at $t$. In this paper, we consider an  $\mathcal{F}$-adapted process $X$ with regulated trajectories of the form $X_{t}=X^{*}_{t}+\sum_{s<t}\Delta^{+}X_{s}, \,\,\,\ t\in[0,T],$
where $X^{*}$ is an $\mathcal{F}$-adapted semimartingale whit c\`{a}dl\`{a}g trajectories and
$\sum_{s<t}|\Delta^{+}X_{s}|<\infty$, $\mathbb{P}-a.s.$  We denote:
\begin{itemize}
\item $\mathcal{T}_{t,T}$ is the set of all  stopping times $\tau$
such that $\mathbb{P}(t\leq \tau \leq T)=1$. More generally, for a
given stopping time  $\nu$ in $\mathcal{T}_{0,T}$, we denote by
   $\mathcal{T}_{\nu,T}$ the set of all  stopping times $\tau$
  such that $\mathbb{P}(\nu \leq \tau \leq T)=1$.
\item $L^{2}(\mathcal{F}_{T})$ is the set of random variables which are
$\mathcal{F}_{T}$-measurable and square-integrable.
\item $\mathcal{M}_{loc}$ is the set of c\`{a}dl\`{a}g local martingales.
\end{itemize}
Now to define the solution of our reflected backward stochastic
differential equation,  let us introduce the following
spaces:
\begin{itemize}
\item $\mathcal{S}^{2}$ is the set of all $\mathcal{F}$-progressively measurable process with regulated trajectories $\phi$ such that:
         $$E\Bigl[\sup_{0\leq t\leq T}|\phi_{t}|^{2}\Bigr]<\infty.$$
\item $\mathbb{M}^{2}$ is the subspace of $\mathcal{M}_{loc}$ of all martingales such that:
         $E([M]_{T})<+\infty.$
\end{itemize}
The random variable $\xi$ is $\mathcal{F}_{T}$-measurable with
values in $\R^{d}$ $(d\geq 1)$ and
$f:\Omega\times[0,T]\times\R^{d}\longrightarrow\R^{d}$ is a random
function measurable with respect to $Prog\times\mathcal{B}(\R^{d})$ where $Prog$ denotes the
$\sigma$-field of progressive subsets
of $\Omega\times[0,T]$. A sequence $\{\tau_{k}\}\subset\mathcal{T}_{0,T}$ is called stationary if
$\forall\omega\in\Omega$,$\exists n\in\N$ $\forall k\geq n$ , $\tau_{k}(\omega)=T.$\\
We will need the following assumptions\\
(H1) There is $\mu\in\R$ such that $|f(t,y)-f(t,y')|\leq\mu|y-y'|$ for all $t\in[0,T]$,  $y,y'\in\R$.\\
(H2) $\xi, \int_{0}^{T}|f(r,0)|dr \in L^{2}$\\
(H3) $[0,T]\ni t\longmapsto f(t,y)\in L^{1}(0,T)$ for every $y\in\R$.\\

\section{Reflected BSDE with one barrier.}

In that follows, we assume that $\xi$ is an $\mathcal{F}$-measurable random variable, $L$  and $U$ are $\mathcal{F}$-adapted optional processes in $\mathcal{S}^{2}$  and $L_{t}\leq U_{t}$, for all $t\leq T$ and $L_{T}\leq\xi\leq U_{T}$. We assume that the lower $L$ obstacle is right upper-semicontinuous (r.u.s.c.) and the  upper obstacle $U$ is right lower-semicontinuous (r.l.s.c.).
\begin{definition}
\label{definition1}
We say that a triple $(Y,M,K)$ of $\mathcal{F}$-progressively measurable processes is a solution of the reflected BSDE with driver $f$, terminal value $\xi$ and lower barrier $L$ $(\underline{RBSDE}(\xi, f, L)$ for short) if
\begin{enumerate}
  \item $Y,K\in\mathcal{S}^{2}$, $M\in\mathcal{M}_{loc}$ with $M_{0}=0$.
  \item  $Y_{t} \geq L_{t}$  for all $t\in[0,T]$ a.s., and
  $\int_{0}^{T}(Y_{s-}-L_{s-})dK_{s}^{*}+\sum_{s<T}(Y_{s}-L_{s})\Delta^{+}K_{s}=0$
  \item $\int_{0}^{T}|f(s,Y_{s})|ds<\infty$ a.s.
  \item $Y_{t}=\xi+\int_{t}^{T}f(s,Y_{s})ds+K_{T}-K_{t}-\int_{t}^{T}dM_{s}$, for all $t\in[0,T]$, a.s.
\end{enumerate}
\end{definition}
\begin{remark}
We note that if $L$ and $K$ are c\`{a}dl\`{a}g, then $(2)$ in Definition 1

 reduces to  $\int_{0}^{T}(Y_{s-}-L_{s-})dK_{s}=0.$
\end{remark}
\begin{definition}
\label{definition2}
We say that a triple $(Y,M,A)$ of $\mathds{F}$-progressively measurable processes is a solution of the reflected BSDE with driver $f$, terminal value $\xi$ and upper barrier $U$ $(\overline{RBSDE}(\xi, f, U)$ for short) if
\begin{enumerate}
  \item  $Y,A\in\mathcal{S}^{2}$, $M\in\mathcal{M}_{loc}$ with $M_{0}=0$.
  \item  $Y_{t}\leq U_{t}$ for all $t\in[0,T]$ a.s., and
  $\int_{0}^{T}(U_{s-}-Y_{s-})dA_{s}^{*}+\sum_{s<T}(U_{s}-Y_{s})\Delta^{+}A_{s}=0$
  \item $\int_{0}^{T}|f(s,Y_{s})|ds<\infty$ a.s.
  \item $Y_{t}=\xi+\int_{t}^{T}f(s,Y_{s})ds-(A_{T}-A_{t})-\int_{t}^{T}dM_{s}$  for all $t\in[0,T]$ a.s.
\end{enumerate}
\end{definition}
\begin{remark}
If $(Y,M,K)\in\mathcal{S}^{2}\times\mathcal{M}_{loc}\times\mathcal{S}^{2}$
satisfies  definition 1  then the process $Y$ has left and right limits. Moreover, the process given by $(Y_{t}+\int_{0}^{t}f(s,Y_{s})ds)_{t\in[0,T]}$ is a strong martingale (\cite[Definition A.1]{Ouknine2015}).
\end{remark}

In the theorem below we recall some results on reflecting BSDEs with one barrier. They will play important role in the proof of our main result.  In the penalization method for reflected BSDEs proposed by Klimsiak et al in \cite{Klimsiaketal2016}, they defined arrays of stopping times $\{\{\sigma_{n,i}\}\}$ exhausting right-side jumps of $L$ inductively as follow:  $\sigma_{1,0}=0$ and
$$\sigma_{1,i}=\inf\{t>\sigma_{1,i-1} : \Delta^{+}L_{s}<-1\}\wedge T, \,\,\,\ i=1,...,k_{1}$$
for some $k_{1}\in\N$. Next for $n\in\N$ and given array $\{\{\sigma_{n,i}\}\}$,  $\sigma_{n+1,0}=0$ and
$$\sigma_{n+1,i}=\inf\{t>\sigma_{n+1,i-1} : \Delta^{+}L_{s}<\frac{-1}{n+1}\}\wedge T$$
for $i=1, ... ,j_{n+1}$ where $j_{n+1}$ is chosen so that $\mathbb{P}(\sigma_{n+1,j_{n+1}} < T)\rightarrow0$ as $n\rightarrow\infty$ and
$$\sigma_{n+1,i}=\sigma_{n+1,j_{n+1}}\vee\sigma_{n,i-j_{n+1}}, \,\,\,\ i= j_{n+1} + 1, ... , k_{n+1}, \,\,\,\ k_{n+1} = j_{n+1} + k_{n}.$$

\begin{theorem}
\label{theorem1}
Assume that $(H1)-(H4)$ are satisfied. Then
\begin{description}
\item[(i)] There exists a unique solution  $(\underline{Y},\underline{M},\underline{K})$ of $\underline{RBSDE}(\xi, f, L)$. Moreover if  $(\underline{Y}^{n},\underline{M}^{n})$, $n\in\N$ are solution of BSDEs of the form
\begin{equation}
\underline{Y}^{n}_{t}=\xi+\int_{t}^{T}f(s,\underline{Y}_{s}^{n})ds-\int_{t}^{T}d\underline{M}^{n}_{s}+n\int_{t}^{T}(\underline{Y}_{s}^{n}-L_{s})^{-}ds+
      \sum_{t\leq\sigma_{n,i}<T}(\underline{Y}^{n}_{\sigma_{n,i}^{+}}-L_{\sigma_{n,i}})^{-}
      \label{equation1}
\end{equation}
      then   $\underline{Y}^{n}_{t}\nearrow \underline{Y}_{t}$, $t\in[0,T]$ $\mathbb{P}$-a.s.
\item[(ii)] There exists a unique solution  $(\overline{Y},\overline{M},\overline{A})$ of $\overline{RBSDE}(\xi, f, U)$. Moreover if  $(\overline{Y}^{n},\overline{M}^{n})$, $n\in\N$ are solution of BSDEs of the form
\begin{equation}
    \overline{Y}^{n}_{t}=\xi+\int_{t}^{T}f(s,\overline{Y}_{s}^{n})ds-\int_{t}^{T}d\overline{M}^{n}_{s}-n\int_{t}^{T}(U_{s}-\overline{Y}_{s}^{n})^{-}ds-
      \sum_{t\leq\sigma_{n,i}<T}(U_{\sigma_{n,i}}-\overline{Y}^{n}_{\sigma_{n,i}^{+}})^{-}
      \label{equation2}
\end{equation}
      then   $\overline{Y}^{n}_{t}\nearrow \overline{Y}_{t}$, $t\in[0,T]$ $\mathbb{P}$-a.s.
\end{description}
\end{theorem}
\begin{proof}
The first part in $(i)$ is proved in  \cite{baadi2018} (see also \cite{Klimsiaketal2016} $(p>1)$ and \cite{Ouknine2015} $(p=2)$  in the case of Brownian filtration and  \cite{baadi2017} for the case of a filtration that supports a Brownian motion  and an independent Poisson random measure).

The second part in $(i)$ is proved for the case of a Brownian filtration in \cite[Theorem 4.1]{Klimsiaketal2016}. To show the results in a general filtration we use the It\^{o} formula for the regulated process (see \cite[Theorem 2.5]{baadi2017} or   \cite[Appendix]{Klimsiaketal2016}) to get this inequality:
\begin{eqnarray*}
  \int_{\sigma}^{\tau}d[M-M^{n}]^{c}_{s} &\leq& |Y_{\tau}-Y_{\tau}^{n}|^{2}
   +2\int_{\sigma}^{\tau}(Y_{s-}-Y_{s-}^{n})(f(s,Y_{s})-f(s,Y^{n}_{s}))ds + 2\int_{\sigma}^{\tau}(Y_{s-}-Y_{s-}^{n})d(K_{s}-K^{n}_{s})^{*}\\ &-& 2\int_{\sigma}^{\tau}(Y_{s-}-Y_{s-}^{n})d(M_{s}-M^{n}_{s})
           - 2\sum_{\sigma\leq s<\tau}(Y_{s}-Y_{s}^{n})\Delta^{+}(Y_{s}-Y_{s}^{n})
\end{eqnarray*}
with $(Y^{n},M^{n})$ defined in $~\eqref{equation1}$,  $\sigma,\tau \in \mathcal{T}_{0,T}$,  $\sigma\leq\tau$, and $K^{n}_{t}=n\int_{0}^{t}(Y_{s}^{n}-L_{s})^{-}ds+
      \sum_{0\leq\sigma_{n,i}\leq t}(Y_{\sigma_{n,i}^{+}}-L_{\sigma_{n,i}})^{-}$.
By the fact that  $\Delta^{+}(Y_{s}-Y_{s}^{n})=-\Delta^{+}(K_{s}-K_{s}^{n})$, we have
\begin{eqnarray*}
  \int_{\sigma}^{\tau}d[M-M^{n}]^{c}_{s}
  &\leq& |Y_{\tau}-Y_{\tau}^{n}|^{2}
   +2\int_{\sigma}^{\tau}(Y_{s-}-Y_{s-}^{n})(f(s,Y_{s})-f(s,Y^{n}_{s}))ds + 2\int_{\sigma}^{\tau}(Y_{s-}-Y_{s-}^{n})d(K_{s}-K^{n}_{s})^{*} \\ &-& 2\int_{\sigma}^{\tau}(Y_{s-}-Y_{s-}^{n})d(M_{s}-M^{n}_{s})
           + 2\sum_{\sigma\leq s<\tau}(Y_{s}-Y_{s}^{n})\Delta^{+}(K_{s}-K_{s}^{n})
\end{eqnarray*}
Then $E\int_{\sigma}^{\tau}d[M-M^{n}]^{c}_{s}\leq E|Y_{\tau}-Y_{\tau}^{n}|^{2}+2E\int_{\sigma}^{\tau}|Y_{s-}-Y_{s-}^{n}|(f(s,Y_{s})-f(s,Y^{n}_{s}))ds$.\\
And with the Theorem 4.1 assumptions in \cite{Klimsiaketal2016} we get the existence of a stationary sequence $\{\tau_{k}\}$ of stopping times such that: $E\int_{0}^{\tau_{k}}d[M-M^{n}]_{s}=E\int_{0}^{\tau_{k}}d[M-M^{n}]^{c}_{s}\rightarrow0$. \\
Therefore to prove that $Y^{n}\nearrow Y, t\in[0,T]$, it suffices to repeat step by step the proof of \cite[Theorem 4.1]{Klimsiaketal2016}.\\
The assertion (ii) follows directly from (i). Indeed, let $(Y, M, A)$ be the  solution for the $\overline{RBSDE}(\xi,f,U)$. Set $f'(s,y)= -f(s,-y)$. Then the process $(-Y,-M,A)$ is the solution of the reflected BSDE $\underline{RBSDE}(-\xi,f',-U)$.

\end{proof}

\section{BSDEs with two reflecting barriers.}

In this section $\xi$, $f$, $L$ and $U$ are as in above. We also suppose that $L_{t}\leq U_{t}$ for $t\in[0,T]$ $\mathbb{P}$-a.s.
\begin{definition}
\label{definition3}
We say that a quadruplet $(Y,M,K,A)$ of $\mathds{F}$-progressively measurable processes is a solution of the reflected BSDE with driver $f$, terminal value $\xi$, lower barrier $L$ and  upper barrier $U$,  $(RBSDE(\xi, f, L, U)$ for short), if
\begin{description}
  \item[(LU1)] $Y,K,A\in\mathcal{S}^{2}$, $M\in\mathcal{M}_{loc}$ with $M_{0}=0$.
  \item[(LU2)]  $L_{t}\leq Y_{t}\leq U_{t}$, $t\in[0,T]$ $\mathbb{P}$-a.s.
  \item[(LU3)] $\int_{0}^{T}(U_{s-}-Y_{s-})dA_{s}^{*}+\sum_{s<T}(U_{s}-Y_{s})\Delta^{+}A_{s}=
      \int_{0}^{T}(Y_{s-}-L_{s-})dK_{s}^{*}+\sum_{s<T}(Y_{s}-L_{s})\Delta^{+}K_{s}=0$, a.s.
  \item[(LU4)] $Y_{t}=\xi+\int_{t}^{T}f(s,Y_{s})ds+(K_{T}-K_{t})-(A_{T}-A_{t})-\int_{t}^{T}dM_{s}$, $t\in[0,T]$ $\mathbb{P}$-a.s.
\end{description}
\end{definition}
\begin{remark}
We note  that, due to equation $(LU4)$, we have $\bigtriangleup^{+}Y_{t}=-\bigtriangleup^{+}(K_{t}-A_{t})$
\end{remark}
We are now going to focus on the uniqueness of the solution of the doubly reflected BSDE associated
with $(f,\xi,L,U)$.  However the first step is to provide a comparison result
between the components $Y$ of two solutions (in Definition 3). Actually we have:
\begin{proposition}
\label{proposition1}
Let $(f,\xi,L,U)$ and $(f',\xi',L',U')$  be two sets of data satisfying
(H1)-(H3). Let $(Y,M,K,A)$ and $(Y',M',K',A')$ be two solutions of the
doubly reflected BSDE associated with $(f,\xi,L,U)$ and $(f',\xi',L',U')$ respectively. Assume that $\xi\leq \xi'$, $L\leq L'$, $U\leq U'$ and $f\leq f'$. Then $ \mathbb{P}-a.s. \,\,\,\ Y_{t}\leq Y'_{t}$.
\end{proposition}

\begin{proof}
 Let $(\tau_{k})_{k\geq0}$ be a non-decreasing  sequence, of stationary type and converges to $T$ such that:
$$\tau_{k}=\inf\{t\geq 0,[M]_{t}+[M']_{t})\geq k \}\wedge T.$$
we have $\mathbb{P}-a.s$.,  $[M]_{T}+[M']_{T}<\infty$. Now, by  It\^{o}-Tanaka's formula  for the regulated process (see \cite[Section 3, page 539]{Lenglart1980})  with $(Y-Y')^{+}$ on $[t\wedge\tau_{k},\tau_{k}]$ we get:
\begin{eqnarray*}
  (Y_{t\wedge\tau_{k}}-Y'_{t\wedge\tau_{k}})^{+} &\leq&
  (Y_{\tau_{k}}-Y'_{\tau_{k}})^{+} +\int_{t\wedge\tau_{k}}^{\tau_{k}}1_{\{Y_{s-}>Y'_{s-}\}}(f(s,Y_{s})-f'(s,Y'_{s}))ds\\ &-& \int_{t\wedge\tau_{k}}^{\tau_{k}}1_{\{Y_{s-}>Y'_{s-}\}}d(M_{s}-M'_{s}) +\int_{t\wedge\tau_{k}}^{\tau_{k}}1_{\{Y_{s-}>Y'_{s-}\}}d(K_{s}-K'_{s}-A_{s}+A'_{s}).
\end{eqnarray*}
From definition of solution we have $\int_{t\wedge\tau_{k}}^{\tau_{k}}1_{\{Y_{s-}>Y'_{s-}\}}d(K_{s}-K'_{s}-A_{s}+A'_{s})\leq 0$, and by using the  Lipschitz condition of $f$, we have
$$
  (Y_{t\wedge\tau_{k}}-Y'_{t\wedge\tau_{k}})^{+}\leq
  (Y_{\tau_{k}}-Y'_{\tau_{k}})^{+} +\mu\int_{t\wedge\tau_{k}}^{\tau_{k}}(Y_{s}-Y'_{s})^{+}ds- \int_{t\wedge\tau_{k}}^{\tau_{k}}1_{\{Y_{s-}>Y'_{s-}\}}d(M_{s}-M'_{s}) .
$$
where $\mu$ the Lipschitz constant of $f$. Therefore taking expectation, the limit as $k\longrightarrow\infty$, we have $E[(Y_{\tau_{k}}-Y'_{\tau_{k}})^{+}]\longrightarrow E[(Y_{T}-Y'_{T})^{+}]=0$  since $L\leq Y\leq U$ and $L\leq Y'\leq U$. And by using Gronwall's Lemma we get $E[(Y_{t}-Y'_{t})^{+}]=0$ for any $t\leq T$, a.s., $Y_{t}\leq Y'_{t}$, which is the desired result.
\end{proof}
\begin{proposition}
The RBSDE$(\xi, f, L, U)$  has at most one solution, i.e., if   $(Y,M,K,A)$  and $(Y',M',K',A')$ are two solutions of RBSDE$(\xi, f, L, U)$, then, $\mathbb{P}-a.s.$, $Y=Y'$, $M=M'$ and $K-A=K'-A'$.
\end{proposition}
\begin{proof}
Let $(Y,M,K,A)$  and $(Y',M',K',A')$ be two solutions of RBSDE$(f,\xi,L,U)$. Then from the comparison result( Proposition 1), we have  $Y_{t}=Y'_{t}$, $t\leq T$, $\mathbb{P}$-a.s. and then   $M=M'$	 and by (\textbf{LU4}) in Definition 3, we get $K-A=K'-A'$.

\end{proof}
\begin{remark}
We have also $K=K'$  and $A = A'$   since $L_{t} < U_{t}$, $\forall t < T$. (see, \cite[Proposition 2.1.]{Asri2011}).
\end{remark}

\section{Local solutions of BSDEs with two optional reflecting barriers}

We are going to construct a solution for the doubly reflected BSDE associated with $(f,\xi,L,U)$  step by step under $(H1)-(H3)$. For this we need to construct a process $Y$ which satisfies
locally the RBSDE$(f,\xi,L,U)$, that is to say, for any stopping time $\tau$, we can find another
appropriate stopping time $\lambda_{\tau}$ such that between $\tau$ and $\lambda_{\tau}$, $Y$ satisfies the doubly reflected BSDE. This local solution will be constructed as a limit of a penalization scheme, which leads to study BSDEs with one reflecting barrier.  Thus our first task is to provide the results we need later on BSDEs with one reflecting barrier.
We first introduce the notion of a local solution of the RBSDE$(f,\xi,L,U)$.
\begin{definition}
 \label{definition4}
Let $\tau$ and $\sigma$ be two stopping times such that $\tau\leq\sigma$ $\mathbb{P}$-a.s..
We say that $(Y_t,M_t,K_t,A_t)_{t\leq T}$ is a local solution on $[\tau,\sigma]$  for the doubly reflected
BSDE associated with two barriers $L$ and $U$, the terminal condition $\xi$ and the
generator $f$ if: $\mathbb{P}$-a.s.,$Y,K,A\in\mathcal{S}^{2}$, $M\in\mathcal{M}_{loc}$, $M_{0}=0$ and
\begin{multline}
\label{eqdefinition4}
\left\{
  \begin{array}{ll}
    Y_{T}=\xi,  \forall t\in[\tau,\sigma]  \\
    Y_{t}=Y_{\sigma}+\int_{t}^{\sigma}f(s,Y_{s})ds+(K_{\sigma}-K_{t})-(A_{\sigma}-A_{t})-\int_{t}^{\sigma}dM_{s}, t\in[\tau,\sigma]\,\,\ \mathbb{P}-a.s.,  \\
 L_{t}\leq Y_{t}\leq U_{t}, \forall t\in[\tau,\sigma], \\
    \int_{\tau}^{\sigma}(U_{s-}-Y_{s-})dA_{s}^{*}+\sum_{\tau\leq s<\sigma}(U_{s}-Y_{s})\Delta^{+}A_{s}=\int_{\tau}^{\sigma}(Y_{s-}-L_{s-})dK_{s}^{*}+\sum_{\tau\leq s<\sigma}(Y_{s}-L_{s})\Delta^{+}K_{s}=0, a.s.
\end{array}
\right.
\end{multline}
\end{definition}
In this section, we are going to show the existence of an appropriate local solution which
later will allow us to construct a global solution for the RBSDE$(f,\xi, L, U)$ with  regulated processes in a  general filtration. But we assume only that $L$ is right upper-semicontinuous (r.u.s.c) and $U$ is right lower-semicontinuous (r.l.s.c).\\
The idea of the proof is the same as in the paper of Hamad\`{e}ne and Hassani \cite{HamadeneHassani2005}, in which the authors  proved the results for the double RBSDE with continuous processes and Brownian filtration.
\subsection{The increasing penalization scheme}
Let us introduce the following increasing penalization scheme. For any $n\geq0$, let $(Y_{t}^{n},M_{t}^{n},A_{t}^{n})$ be the triple of $\mathcal{F}_{t}$-adapted processes with values in $\R\times\R^{d}\times\R$, unique solution of the RBSDE$(f(s,y)+n(y-L_{s})^{-}+\sum_{\sigma_{n,i}<t}(y_{\sigma_{n,i}^{+}}-L_{\sigma_{n,i}})^{-},\xi, U)$ such that:
$Y^{n}$,$A^{n}\in\mathcal{S}^{2}$, $M^{n}\in\mathcal{M}_{loc}$ with $M_{0}^{n}=0$ and
\begin{multline}
\label{sheme1}
\left\{
  \begin{array}{ll}
    Y^{n}_{t}=\xi+\int_{t}^{T}(f(s,Y^{n}_{s})+n(Y^{n}_{s}-L_{s})^{-})ds+\sum_{t\leq\sigma_{n,i}<T}(Y^{n}_{\sigma_{n,i}^{+}}-L_{\sigma_{n,i}})^{-} -A^{n}_{T}+A^{n}_{t}-\int_{t}^{T}dM^{n}_{s}, \,\,\ \text{a.s.}, \\
   \forall t\in[0,T], \,\,\ Y_{t}^{n}\leq U_{t} \,\,\  \text{and} \,\,\ \int_{0}^{T}(U_{s-}-Y^{n}_{s-})dA^{n,*}_{s}+\sum_{s<T}(U_{s}-Y^{n}_{s})\Delta^{+}A^{n}_{s}=0, \,\,\ \text{a.s.}
\end{array}
\right.
\end{multline}
We set $f_{n}(t,y)=f(t,y)+n(y-L_{t})^{-}+\sum_{\sigma_{n,i}<t}(y_{\sigma_{n,i}^{+}}-L_{\sigma_{n,i}})^{-}.$\\
By Theorem 1 there exist a unique solution $(Y_{t}^{n},M_{t}^{n},A_{t}^{n})$ of RBSDE$(f_{n}(s,y),\xi, U)$. We have $f_{n}(s,y)\leq f_{n+1}(s,y)$ which implies from the
comparison result that for any $n\geq 0$, we have $Y^{n}\leq Y^{n+1} \leq U$. And by consequence there exist  $Y=(Y_{t})_{t\leq T}$ such that
$(Y^{n}_{t})_{t\leq T}$ converges increasingly to $(Y_{t})_{t\leq T}$ and for any $t\leq T$, $Y_{t}\leq U_{t}$.
Besides for a stopping time $\tau$ let us set: $\delta_{\tau}^{n}=\inf\{s\geq \tau, Y_{s}^{n}=U_{s}\}\wedge T.$\\
Since  $Y^{n}\leq Y^{n+1}$ then the sequence $(\delta_{\tau}^{n})_{n\geq 0}$ is decreasing and converges to $\delta_{\tau}$.  Let us now  focus on some
properties of $Y$ and especially show that $Y\geq L$.
\begin{proposition}
\label{proposition3}
The following properties are fulfilled $\mathbb{P}$-a.s.:
\begin{description}
  \item[(i)] $Y_{\delta_{\tau}}1_{[\delta_{\tau}<T]}=U_{\delta_{\tau}}1_{[\delta_{\tau}<T]}.$
  \item[(ii)] $\forall t\leq T, Y_{t}\geq L_{t}.$
\end{description}
\end{proposition}
\begin{proof}
We begin with the proof of (i). For $n\geq0$ and $t\leq T$  the process $A^{n}$ does not increase before $Y^{n}$ reaches the barrier $U$, then for any $t\in[\tau,\delta_{\tau}^{n}]$, we have, $A^{n}_{t}-A^{n}_{\tau}=0$ and then
\begin{equation}
Y^{n}_{t}=Y^{n}_{\delta_{\tau}^{n}}+\int_{t}^{\delta_{\tau}^{n}}(f(s,Y^{n}_{s})+n(Y^{n}_{s}-L_{s})^{-})ds+\sum_{t\leq\sigma_{n,i}<\delta_{\tau}^{n}}(Y^{n}_{\sigma_{n,i}^{+}}-L_{\sigma_{n,i}})^{-} -\int_{t}^{\delta_{\tau}^{n}}dM^{n}_{s}.
      \label{equation3}
\end{equation}
For $n\geq0$, writing $~\eqref{equation3}$ between $\delta_{\tau}$ and $\delta_{\tau}^{n}$ ($\delta_{\tau}^{n}\searrow \delta_{\tau}$) yields:
\begin{equation}
Y^{n}_{\delta_{\tau}}=Y^{n}_{\delta_{\tau}^{n}}+\int_{\delta_{\tau}}^{\delta_{\tau}^{n}}(f(s,Y^{n}_{s})+n(Y^{n}_{s}-L_{s})^{-})ds+
\sum_{\delta_{\tau}\leq\sigma_{n,i}<\delta_{\tau}^{n}}(Y^{n}_{\sigma_{n,i}^{+}}-L_{\sigma_{n,i}})^{-} -\int_{\delta_{\tau}}^{\delta_{\tau}^{n}}dM^{n}_{s}
      \label{equation4}
\end{equation}
and then $Y_{\delta_{\tau}}\geq U_{\delta_{\tau}^{n}}1_{[\delta_{\tau}^{n}<T]}+\xi1_{[\delta_{\tau}^{n}=T]}+\int_{\delta_{\tau}}^{\delta_{\tau}^{n}}f(s,Y^{n}_{s})ds -\int_{\delta_{\tau}}^{\delta_{\tau}^{n}}dM^{n}_{s}$, which implies that
\begin{equation}
1_{[\delta_{\tau}<T]}Y_{\delta_{\tau}}\geq 1_{[\delta_{\tau}<T]}(U_{\delta_{\tau}^{n}}1_{[\delta_{\tau}^{n}<T]}+\xi1_{[\delta_{\tau}^{n}=T]})+1_{[\delta_{\tau}<T]}\int_{\delta_{\tau}}^{\delta_{\tau}^{n}}f(s,Y^{n}_{s})ds -\int_{\delta_{\tau}}^{\delta_{\tau}^{n}}1_{[\delta_{\tau}<T]}dM^{n}_{s}
      \label{equation5}
\end{equation}
By (H1) we have: $|f(s,Y^{n}_{s})|\leq |f(s,0)|+\mu|Y^{n}_{s}|$ where $\mu$ is a constant. We have also  $Y^{0}\leq Y^{n}\leq U$ which implies that $E[\int_{\delta_{\tau}}^{\delta_{\tau}^{n}}|Y^{n}_{s}|ds]$ converges to $0$. Consequently, $\limsup_{n\longrightarrow\infty}E\Bigl[\int_{\delta_{\tau}}^{\delta_{\tau}^{n}}|f(s,Y^{n}_{s})|ds\Bigr]=0.$\\
Using now  inequality $~\eqref{equation5}$ and taking expectation in both hand-sides then
the limit as $n$ goes to infinity to obtain:
$E\Bigl[1_{[\delta_{\tau}<T]}Y_{\delta_{\tau}}\Bigr]\geq E\Bigl[1_{[\delta_{\tau}<T]}U_{\delta_{\tau+}}\Bigr]\geq E\Bigl[1_{[\delta_{\tau}<T]}U_{\delta_{\tau}}\Bigr],$
since $U$ is optional r.l.s.c. process. By  $Y\leq U$, we have the desired result.

We now prove (ii). For any $n\geq 0$ and any stopping time $\tau\leq T$, the following property holds true:
$$Y^{n}_{\tau}=E\Bigl[\int_{\tau}^{\delta_{\tau}^{n}}(f(s,Y^{n}_{s})+n(Y^{n}_{s}-L_{s})^{-})ds+
\sum_{\tau\leq\sigma_{n,i}<\delta_{\tau}^{n}}(Y^{n}_{\sigma_{n,i}^{+}}-L_{\sigma_{n,i}})^{-}+
U_{\delta_{\tau}^{n}}1_{[\delta_{\tau}^{n}<T]}+\xi1_{[\delta_{\tau}^{n}=T]}/\mathcal{F}_{\tau}\Bigr],$$
since the process $A^{n}$ does not increase before $Y^{n}$ reaches the barrier $U$  by definition of  $\delta_{\tau}^{n}$. From last equality  we have
\begin{multline}
E\Bigl[\int_{\tau}^{\delta_{\tau}^{n}}(Y^{n}_{s}-L_{s})^{-}ds+
\frac{1}{n}\sum_{\tau\leq\sigma_{n,i}<\delta_{\tau}^{n}}(Y^{n}_{\sigma_{n,i}^{+}}-L_{\sigma_{n,i}})^{-}\Bigr]
\leq \frac{1}{n} E\Bigl[\int_{\tau}^{\delta_{\tau}^{n}}|f(s,Y^{n}_{s})|ds +|Y^{n}_{\tau}|+|U_{\delta_{\tau}^{n}}|1_{[\delta_{\tau}^{n}<T]}+|\xi|1_{[\delta_{\tau}^{n}=T]}\Bigr]
\label{equation6}
\end{multline}
By (H3) we have  $E[\int_{0}^{T}|f(s,Y^{n}_{s})|ds]<\infty$ when $n$ goes to infinity, and by Fatou's lemma
 we deduce from $~\eqref{equation6}$ that:
$E\Bigl[\int_{0}^{T}\liminf_{n\longrightarrow\infty}1_{[\tau,\delta_{\tau}^{n}]}(Y_{s}-L_{s})^{-}ds\Bigr]\leq
\liminf_{n\longrightarrow\infty}E\Bigl[\int_{\tau}^{\delta^{n}_{\tau}}(Y_{s}-L_{s})^{-}ds\Bigr]=0,$ then
\begin{equation}
\int_{\tau}^{\delta_{\tau}}(Y_{s}-L_{s})^{-}ds=0 \,\,\,\ a.s.
      \label{equation7}
\end{equation}

Since  $Y^{n}\leq Y^{n+1}$, note that if  $L$ is a c\`{a}dl\`{a}g process the limit $Y$ of $\{Y^{n}\}$ is c\`{a}dl\`{a}g  (\cite[Theorem 3.1]{Essaky2008}  and  \cite[Lemma 2.2]{Peng1999}) on $[\tau,\delta_{\tau}]$. But in our case $Y$ need not to be c\`{a}dl\`{a}g. Henceforth from $~\eqref{equation7}$ we obtain that $Y_{\tau}\geq L_{\tau}$ on the set $[\tau,\delta_{\tau}]$. If  $\tau=\delta_{\tau}<T$ we have $Y_{\tau}=U_{\tau}\geq L_{\tau}$ and if $\tau=\delta_{\tau}=T$ we have $Y_{\tau}=\xi\geq L_{\tau}$. By consequence for all $\tau$, $Y_{\tau}\geq L_{\tau}$. As
the barriers $L$ and $U$ are optional, using the optional section theorem \cite[Proposition A.4]{Ouknine2015} we have $\mathbb{P}$-a.s., $Y\geq L$. The proof is complete.
\end{proof}
Next, we show  the existence of the local solution of the reflected BSDE$(f,\xi,L,U)$ on $[\tau,\delta_{\tau}]$.
\begin{proposition}
\label{proposition4}
There exists two measurable processes $(\bar{K}^{\tau}_{t})_{t\leq T}$ and $(\bar{M}^{\tau}_{t})_{t\leq T}$ such that $(Y_{t},\bar{M}^{\tau}_{t}, \bar{K}^{\tau}_{t}, 0)_{t\leq T}$ is a local solution of RBSDE  in Definition~\ref{definition3}  on $[\tau,\delta_{\tau}]$. That is: $\bar{K}^{\tau}\in\mathcal{S}^{2}$, $\bar{M}^{\tau}\in\mathcal{M}_{loc}$, $\bar{M}^{\tau}_{0}=0$ and
\begin{multline}
\label{eqproposition4}
$$\left\{
  \begin{array}{ll}
    Y_{t}=Y_{\delta_{\tau}}+\int_{t}^{\delta_{\tau}}f(s,Y_{s})ds+\bar{K}^{\tau}_{\delta_{\tau}}-\bar{K}^{\tau}_{t}-\int_{t}^{\delta_{\tau}}d\bar{M}^{\tau}_{s},a.s. \,\,\,\ \text{and} \,\,\,\ Y_{T}=\xi,\\
 \forall t\in[\tau,\delta_{\tau}], L_{t}\leq Y_{t}\leq U_{t} \,\,\,\ \text{and} \,\,\,\
     \int_{\tau}^{\delta_{\tau}}(Y_{s-}-L_{s-})d\bar{K}^{\tau,*}_{s}+\sum_{\tau\leq s<\delta_{\tau}}(Y_{s-}-L_{s-})\Delta^{+}\bar{K}^{\tau}_{s}=0, a.s..
\end{array}
\right.$$
\end{multline}

\end{proposition}
\begin{proof}
For any $n\geq 0$ and $t\in[\tau,\delta_{\tau}]$ and
since the process $A^{n}$ moves only when $Y^{n}$ reaches the barrier $U$ (then $A^{n}_{\tau}=A^{n}_{\delta_{\tau}}$),  we have
$$Y_{t}^{n}=Y^{n}_{\delta_{\tau}}+\int_{t}^{\delta_{\tau}}f(s,Y^{n}_{s})ds+
n\int_{t}^{\delta_{\tau}}(Y^{n}_{s}-L_{s})^{-}ds+
\sum_{t\leq\sigma_{n,i}<\delta_{\tau}}(Y^{n}_{\sigma_{n,i}^{+}}-L_{\sigma_{n,i}})^{-}
-\int_{t}^{\delta_{\tau}}dM^{n}_{s}.$$
On the other hand, for $n\geq 0$, let $(\bar{Y}^{n},\bar{M}^{n})_{t\leq \delta_{\tau}}$ be the unique solution of the BSDE associated with the coefficient $f(t,y_{t})+n(y_{t}-L_{t})^{-}+
\sum_{\sigma_{n,i}<t}(y_{\sigma_{n,i}^{+}}-L_{\sigma_{n,i}})^{-}$, the terminal value $Y^{n}_{\delta_{\tau}}$ and a
bounded terminal time $\delta_{\tau}$, that is,
\begin{multline}
\left\{
  \begin{array}{ll}
     E(\sup_{s\leq\delta_{\tau}}\mid\bar{Y}^{n}_{s}\mid^{2}+ [\bar{M}^{n}]_{\delta_{\tau}})<\infty \\
    \bar{Y}^{n}_{t}=Y^{n}_{\delta_{\tau}}+\int_{t}^{\delta_{\tau}}f(s,\bar{Y}^{n}_{s})ds+\int_{t}^{\delta_{\tau}}n(\bar{Y}^{n}_{s}-L_{s})^{-}ds+
\sum_{t\leq\sigma_{n,i}<\delta_{\tau}}(\bar{Y}^{n}_{\sigma_{n,i}^{+}}-L_{\sigma_{n,i}})^{-}-\int_{t}^{\delta_{\tau}}d\bar{M}^{n}_{s}.
\end{array}
\right.
\label{equation8}
\end{multline}
The proof of existence and uniqueness is obtained  by the same arguments such that in \cite[3.2. Proposition]{HamadeneHassani2005} or \cite[Proposition. 4]{Hassairi2016} since $\delta_{\tau}$ is bounded. We have  $(Y^{n}_{\delta_{\tau}})_{n\geq 0}\nearrow Y_{\delta_{\tau}}\leq U_{\delta_{\tau}}$, hence
from the Lebesgue dominated convergence theorem we get $E(\mid Y^{n}_{\delta_{\tau}}-Y_{\delta_{\tau}}\mid)=0$ as $n\rightarrow\infty$. Therefore the sequence of processes $((\bar{Y}^{n}_{t},\bar{M}^{n}_{t},\int_{0}^{t}n(\bar{Y}^{n}_{s}-L_{s})^{-}ds+
\sum_{\sigma_{n,i}<t}(\bar{Y}^{n}_{\sigma_{n,i}^{+}}-L_{\sigma_{n,i}})^{-})_{t\leq \delta_{\tau}})_{n\geq 0}$ converges in $\mathcal{S}^{2}_{\delta_{\tau}}\times\mathbb{M}^{2}_{\delta_{\tau}}\times\mathcal{S}^{2}_{\delta_{\tau}}$ ($\mathcal{S}^{2}_{\delta_{\tau}}$ and $\mathbb{M}^{2}_{\delta_{\tau}}$ are the same as $\mathcal{S}^{2}$ and $\mathbb{M}^{2}$ except for that $T$ is replaced by the
stopping time $\delta_{\tau}$) to  $(\hat{Y}_{t},\hat{M}_{t},\hat{K}_{t})_{t\leq \delta_{\tau}}$ such that:
$$\left\{
  \begin{array}{ll}
     E(\sup_{s\leq\delta_{\tau}}\mid\hat{Y}_{s}\mid^{2}+ [\hat{M}]_{\delta_{\tau}})<\infty, \hat{K}_{s}\in\mathcal{S}^{2}_{\delta_{\tau}}  \,\,\,\ \text{and} \,\,\,\  \hat{K}_{0}=0. \\
    \hat{Y}_{t}=Y_{\delta_{\tau}}+\int_{t}^{\delta_{\tau}}f(s,\hat{Y}_{s})ds+(\hat{K}_{\delta_{\tau}}-\hat{K}_{t})-\int_{t}^{\delta_{\tau}}d\hat{M}_{s}, \forall t\leq \delta_{\tau},  \\
\hat{Y}_{t}\geq L_{t} \,\,\,\ \text{and} \,\,\,\  \int_{\tau}^{\delta_{\tau}}(\hat{Y}_{s-}-L_{s-})d\hat{K}_{s}^{*}+\sum_{\tau\leq s<\delta_{\tau}}(\hat{Y}_{s-}-L_{s-})\Delta^{+}\hat{K}_{s}=0.
\end{array}
\right.$$
Now by $~\eqref{sheme1}$, $~\eqref{equation8}$ and uniqueness of the solution  on $[\tau,\delta_{\tau}]$ implies that for any $t\in[\tau, \delta_{\tau}]$, $Y^{n}_{t}=\bar{Y}^{n}_{t}$ and $M^{n}_{t}=\bar{M}^{n}_{t}$.
Therefore $Y_{t}=\hat{Y}_{t}$ for any  $t\in[\tau, \delta_{\tau}]$,
$$\left\{
  \begin{array}{ll}
     E(\sup_{s\leq\delta_{\tau}}\mid Y_{s}\mid^{2}+ [\hat{M}]_{\delta_{\tau}})<\infty; \hat{K}_{s}\in\mathcal{S}^{2}_{\delta_{\tau}} \,\,\,\ \text{and} \,\,\,\ \hat{K}_{0}=0 \\
    Y_{t}=Y_{\delta_{\tau}}+\int_{t}^{\delta_{\tau}}f(s,Y_{s})ds+(\hat{K}_{\delta_{\tau}}-\hat{K}_{t})-\int_{t}^{\delta_{\tau}}d\hat{M}_{s}, \forall t\leq \delta_{\tau}, \\
\forall t\leq \delta_{\tau} L_{t}\leq Y_{t}\leq U_{t}  \,\,\,\ \text{and} \,\,\,\ \int_{\tau}^{\delta_{\tau}}(Y_{s-}-L_{s-})d\hat{K}_{s}^{*}+\sum_{\tau\leq s<\delta_{\tau}}(Y_{s-}-L_{s-})\Delta^{+}\hat{K}_{s}=0.
\end{array}
\right.$$
For any $t\leq T$, let us set $\bar{K}^{\tau}_{t}=(\hat{K}_{t\wedge\delta_{\tau}}-\hat{K}_{\tau})1_{[t\geq\tau]}$ and $\bar{M}^{\tau}_{t}=\hat{M}_{t}1_{[\tau\leq t\leq\delta_{\tau}]}$ (see Remark~\ref{remark8}), we
deduce that  $(Y_{t},\bar{M}^{\tau}_{t}, \bar{K}^{\tau}_{t}, 0)_{t\leq T}$ is a local solution of RBSDE  in Definition~\ref{definition3} on $[\tau,\delta_{\tau}]$.
\end{proof}

\subsection{The decreasing penalization scheme}

We now consider the following decreasing penalization scheme for any $n\geq 0$:
\begin{multline}
\label{sheme2}
\left\{
  \begin{array}{ll}
    \tilde{Y}^{n},\tilde{K}^{n}\in\mathcal{S}^{2}, \tilde{M}^{n}\in\mathcal{M}_{loc} \,\,\ with \,\,\ \tilde{M}_{0}^{n}=0.\\
    \tilde{Y}^{n}_{t}=\xi+\int_{t}^{T}(f(s,\tilde{Y}^{n}_{s})-n(U_{s}-\tilde{Y}^{n}_{s})^{-})ds-\sum_{t\leq\sigma_{n,i}<T}(U_{\sigma_{n,i}}-\tilde{Y}^{n}_{\sigma_{n,i}^{+}})^{-} +\tilde{K}^{n}_{T}-\tilde{K}^{n}_{t}-\int_{t}^{T}d\tilde{M}^{n}_{s}, \,\,\ a.s.,\\
    \forall t\in[0,T],\,\,\ \tilde{Y}_{t}^{n}\geq L_{t}, \,\,\ and \,\,\ \int_{0}^{T}(\tilde{Y}^{n}_{s-}-L_{s-})d\tilde{K}^{n*}_{s}+\sum_{s<T}(\tilde{Y}^{n}_{s}-L_{s})\Delta^{+}\tilde{K}^{n}_{s}=0, \,\,\ a.s.
\end{array}
\right.
\end{multline}
First we note that the existence of the triple $(\tilde{Y}^{n},\tilde{M}^{n},\tilde{K}^{n})$  is  due to \cite[Theorem 4.1]{Klimsiaketal2016}  and the following remark.
\begin{remark}
\label{remark5}
A triple $(Y,M,K)$ is a solution for the BSDE with a lower reflecting
barrier associated with $(f,\xi,L)$ iff $(-Y,-M,K)$ is a solution of the BSDE with an
upper reflecting barrier associated with $(-f(t,\omega,-y),-\xi,-L)$.
\end{remark}
For any stopping time $\tau\leq T$ and any $n\geq0$, let us set $\theta_{\tau}^{n}=\inf\{s\geq\tau, \tilde{Y}^{n}_{s}=L_{s}\}\wedge T.$ By Proposition 1,  we have $\tilde{Y}^{n} \geq \tilde{Y}^{n+1} \geq L$  then the sequence
$(\tilde{Y}^{n})_{n\geq0}$ converges to  $\tilde{Y}$ and $(\theta_{\tau}^{n})_{n\geq0}$ is decreasing
and converges to another stopping time $\theta_{\tau}= \lim_{n\longrightarrow\infty}\theta_{\tau}^{n}$. Using the same arguments in Propositions ~\ref{proposition3} and ~\ref{proposition4}  and by Remark~\ref{remark5}, we get:
\begin{proposition}
\label{proposition5}
The following properties hold true $\mathbb{P}$-a.s.:
\begin{description}
  \item[(i)] $\tilde{Y}_{\theta_{\tau}}1_{[\theta_{\tau}<T]}=L_{\theta_{\tau}}1_{[\theta_{\tau}<T]}.$
  \item[(ii)] $\forall t\leq T, \tilde{Y}_{t}\leq U_{t}.$

 \item[(iii)] There exists two measurable processes $(\tilde{A}^{\tau}_{t})_{t\leq T}$ and $(\tilde{M}^{\tau}_{t})_{t\leq T}$ such that $(\tilde{Y}_{t},\tilde{M}^{\tau}_{t},0, \tilde{A}^{\tau}_{t})_{t\leq T}$ is a local solution of RBSDE  in Definition~\ref{definition3} on $[\tau,\theta_{\tau}]$. That is: $\tilde{A}^{\tau}\in\mathcal{S}^{2}$, $\tilde{M}^{\tau}\in\mathcal{M}_{loc}$, $\tilde{M}^{\tau}_{0}=0$ and  \\
\begin{multline}
\label{eqproposition5}
\left\{
  \begin{array}{ll}
    Y_{T}=\xi \,\,\,\,\ and \,\,\,\,\ \tilde{Y}_{t}=\tilde{Y}_{\theta_{\tau}}+\int_{t}^{\theta_{\tau}}f(s,\tilde{Y}_{s})ds-\tilde{A}^{\tau}_{\theta_{\tau}}+\tilde{A}^{\tau}_{t}-\int_{t}^{\theta_{\tau}}d\tilde{M}^{\tau}_{s}, \,\,\,\ \mathbb{P}-a.s.\\
 \forall t\in[\tau,\theta_{\tau}], L_{t}\leq \tilde{Y}_{t}\leq U_{t} \,\,\,\ and \,\,\,\
     \int_{\tau}^{\theta_{\tau}}(U_{s-}-\tilde{Y}_{s-})d\tilde{A}^{\tau *}_{s}+\sum_{\tau\leq s<\theta_{\tau}}(U_{s-}-\tilde{Y}_{s-})\triangle^{+}\tilde{A}^{\tau}_{s}=0, \,\,\,\ a.s..
\end{array}
\right.
\end{multline}
\end{description}
\end{proposition}
\subsection{Existence of the local solution}
Recall that $Y$ (resp. $\tilde{Y}$) is the limit of the increasing (resp. decreasing) approximating scheme. We are going to show that  the processes $Y$ and  $\tilde{Y}$ are undistinguishable.
\begin{proposition}
\label{proposition6}
$\mathbb{P}$-a.s., for any $t\leq T$ , $Y_{t}=\tilde{Y}_{t}$.
\end{proposition}
\begin{proof}
We prove the equality in two steps, first we show that  $Y\leq\tilde{Y}$, and second we show the other inequality. For that, let $J^{0}(Y^{n}-\tilde{Y}^{m})$ denote the local time of  $Y^{n}-\tilde{Y}^{m}$ at $0$. For any $t\leq T$ and any $n,m\geq0$, by the It\^{o}-Tanaka formula for regulated processes (see ~\cite[Section 3, page 539]{Lenglart1980}) applied to  $(Y^{n}-\tilde{Y}^{m})^{+}$ we have

\begin{multline}
(Y^{n}_{t}-\tilde{Y}^{m}_{t})^{+}\leq (Y^{n}_{T}-\tilde{Y}^{m}_{T})^{+}+ \int_{t}^{T}1_{\{Y^{n}_{s-}>\tilde{Y}^{m}_{s-}\}}(f(s,Y^{n}_{s})-f(s,\tilde{Y}^{m}_{s}))ds+\\
\int_{t}^{T}1_{\{Y^{n}_{s-}>\tilde{Y}^{m}_{s-}\}}d(K^{n}_{s}-A^{n}_{s}-\tilde{K}^{m}_{s}+\tilde{A}^{m}_{s})-
\int_{t}^{T}1_{\{Y^{n}_{s-}>\tilde{Y}^{m}_{s-}\}}d(M^{n}_{s}-\tilde{M}^{m}_{s})
\label{TNAKA}
\end{multline}
As in the proof of the comparison result (see Proposition ~\ref{proposition1}) we chow that $Y^{n}_{t}\leq \tilde{Y}^{m}_{t}$ and  we get $Y_{t}\leq \tilde{Y}_{t}$, for any $t\leq T$.\\
Now we prove that $Y_{t}\geq \tilde{Y}_{t}, \forall t\leq T$. Let $\tau$ be a stopping time and let $\delta_{\tau}$ and $\theta_{\tau}$ be the stopping times introduced in Proposition ~\ref{proposition3} and ~\ref{proposition5} respectively. We have:
\begin{multline}
Y_{\delta_{\tau}\wedge\theta_{\tau}}=Y_{\delta_{\tau}}1_{[\delta_{\tau}\leq\theta_{\tau}<T]}
+Y_{\theta_{\tau}}1_{[\delta_{\tau}>\theta_{\tau}]}
+Y_{\delta_{\tau}}1_{[\delta_{\tau}\leq\theta_{\tau}=T]}
\geq L_{\delta_{\tau}}1_{[\delta_{\tau}\leq\theta_{\tau}<T]}
+U_{\theta_{\tau}}1_{[\delta_{\tau}>\theta_{\tau}]}
+Y_{\delta_{\tau}}1_{[\delta_{\tau}\leq\theta_{\tau}=T]}
\\
\geq L_{\delta_{\tau}}1_{[\delta_{\tau}\leq\theta_{\tau}<T]}
+\tilde{Y}_{\theta_{\tau}}1_{[\delta_{\tau}>\theta_{\tau}]}
+\tilde{Y}_{\delta_{\tau}}1_{[\delta_{\tau}\leq\theta_{\tau}=T]}=\tilde{Y}_{\delta_{\tau}\wedge\theta_{\tau}}
\label{equation9}
\end{multline}
since $Y\geq L$, $\tilde{Y}_{\theta_{\tau}}1_{[\delta_{\tau}>\theta_{\tau}]}=U_{\theta_{\tau}}1_{[\delta_{\tau}>\theta_{\tau}]}$ (Proposition~\ref{proposition5}) and
\begin{eqnarray*}
Y_{\delta_{\tau}}1_{[\delta_{\tau}\leq\theta_{\tau}=T]} = Y_{\delta_{\tau}}1_{[\delta_{\tau}\leq\theta_{\tau}=T]\cap[\delta_{\tau}<T]}
+Y_{\delta_{\tau}}1_{[\delta_{\tau}\leq\theta_{\tau}=T]\cap[\delta_{\tau}=T]}\\ \geq L_{\delta_{\tau}}1_{[\delta_{\tau}\leq\theta_{\tau}=T]\cap[\delta_{\tau}<T]}
+\xi1_{[\delta_{\tau}\leq\theta_{\tau}=T]\cap[\delta_{\tau}=T]}=\tilde{Y}_{\delta_{\tau}}1_{[\delta_{\tau}\leq\theta_{\tau}=T]}.
\end{eqnarray*}
Since $Y$ and $\tilde{Y}$ satisfy the BSDEs ~\eqref{eqproposition4} and ~\eqref{eqproposition5}  respectively between $\tau$ and $\delta_{\tau}\wedge\theta_{\tau}$, then using comparison result of solution (Proposition~\ref{proposition1}) of BSDEs with $(1_{[\tau<\delta_{\tau}\wedge\theta_{\tau}]}Y_{\tau})_{t\in[\tau,\delta_{\tau}\wedge\theta_{\tau}]}$ and $(1_{[\tau<\delta_{\tau}\wedge\theta_{\tau}]}\tilde{Y}_{\tau})_{t\in[\tau,\delta_{\tau}\wedge\theta_{\tau}]}$, we get $1_{[\tau<\delta_{\tau}\wedge\theta_{\tau}]}Y_{\tau}\geq1_{[\tau<\delta_{\tau}\wedge\theta_{\tau}]}\tilde{Y}_{\tau}.$
On the other hand from ~\eqref{equation9}, we have
$1_{[\tau=\delta_{\tau}\wedge\theta_{\tau}]}Y_{\tau}\geq1_{[\tau=\delta_{\tau}\wedge\theta_{\tau}]}\tilde{Y}_{\tau}$
which implies that $Y_{\tau}\geq \tilde{Y}_{\tau}$. As $\tau$ is an arbitrary stopping time and $Y$ and $\tilde{Y}$ are optional processes then $\mathbb{P}$-a.s. $Y\geq \tilde{Y}$ by \cite[Proposition 2.4.]{baadi2017}. We conclude that $Y=\tilde{Y}$ $\mathbb{P}$-a.s..
\end{proof}
As a consequence of   the result obtained in Propositions~\ref{proposition4}, ~\ref{proposition5} and ~\ref{proposition6} we have:
\begin{theorem}
\label{theorem2}
There exists a unique  measurable process ($Y_{t})_{t\leq T}$
such that:
\begin{description}
  \item[i] $\forall t \leq T$, $L_{t}\leq Y_{t}\leq U_{t}$ and $Y_{T}=\xi$  \\
  \item[ii] for any stopping time $\tau$, there exist another stopping time $\lambda_{\tau}\geq\tau$ $\mathbb{P}$-a.s., and a triple of measurable processes $(M^{\tau}_{t},K^{\tau}_{t},A^{\tau}_{t})_{t\leq T}$ such that on $[\tau,\lambda_{\tau}]$ the process $(Y_{t},M^{\tau}_{t},K^{\tau}_{t},A^{\tau}_{t})_{t\leq T}$ is a local solution for the reflected BSDE associated  $(f,\xi,L,U)$, ( ~\eqref{eqdefinition4} in Definition~\ref{definition4}).
\end{description}
\end{theorem}
\begin{remark}
If we set $\nu_{\tau}=\inf\{t\geq\tau, Y_{t}=U_{t}\}\wedge T$, $\sigma_{\tau}=\inf\{t\geq\tau, Y_{t}=L_{t}\}\wedge T$,   when $\nu_{\tau}\vee\sigma_{\tau}\leq\lambda_{\tau}$, that is $Y$ reaches $L$ and $U$ between the times $\tau$ and $\lambda_{\tau}$ when $\lambda_{\tau}<T$.
\end{remark}
\begin{proof} By  Proposition~\ref{proposition4} and  Proposition~\ref{proposition5} we have the first point  (i). Let $\tau$ be a fixed stopping time and let $(Y,\lambda_{\tau},M^{\tau}, K^{\tau}, A^{\tau})$ and $(Y',\lambda'_{\tau},M^{'\tau}, K^{'\tau}, A^{'\tau})$ two solutions for the reflected BSDE associated  $(f,\xi,L,U)$. Then by the same argument of \cite[Theorem 3.2.]{Hassairi2016}, we can prove that $Y=Y'$.

Let $(Y_{t},\bar{M}^{\tau}_{t}, \bar{K}^{\tau}_{t}, 0)_{t\leq T}$ (resp. $(Y_{t},\tilde{M}^{\tau}_{t},0, \tilde{A}^{\tau}_{t})_{t\leq T}$)  be a local solution of reflected BSDE in definition~\ref{definition4} on $[\tau,\delta_{\tau} ]$ (resp. on $[\delta_{\tau},\lambda_{\tau}]$) which exists according to Proposition~\ref{proposition4}  (resp. Proposition~\ref{proposition5}), where $\lambda_{\tau}$  is a stopping time such that $\tau\leq\lambda_{\tau}\leq T$. Now for $t\leq T$, let $M_{t}=\bar{M}_{t}1_{[\tau\leq t\leq\delta_{\tau}]}+\tilde{M}_{t}1_{[\delta_{\tau}\leq t\leq\lambda_{\tau}]}$ (see Remark~\ref{remark8}), $K_{t}=\bar{K}_{t\wedge\delta_{\tau}}$  and $A_{t}=\tilde{A}_{t\wedge\delta_{\tau}}1_{[t\geq \delta_{\tau}]}$.
For any $t\in[\tau,\lambda_{\tau}]$ we have,
\begin{multline}
\label{eqtheorem2}
\left\{
  \begin{array}{ll}
    Y_{t}=Y_{\lambda_{\tau}}+\int_{t}^{\lambda_{\tau}}f(s,Y_{s})ds+K^{\tau}_{\lambda_{\tau}}-K^{\tau}_{t}-(A^{\tau}_{\lambda_{\tau}}-A^{\tau}_{t})-\int_{t}^{\lambda_{\tau}}dM^{\tau}_{s}, \,\,\,\  \mathbb{P}-a.s.\\
 \forall t\in[\tau,\lambda_{\tau}], L_{t}\leq Y_{t}\leq U_{t} \,\,\,\ \text{a.s.}, \,\,\,\ \text{and}  \\ \int_{\tau}^{\lambda_{\tau}}(Y_{s-}-L_{s-})dK^{\tau *}_{s}+\sum_{\tau\leq s<\lambda_{\tau}}(Y_{s-}-L_{s-})\Delta^{+}K^{\tau}_{s}=0 \,\,\,\ \text{a.s.}, \,\,\,\ \text{and}\\ \int_{\tau}^{\lambda_{\tau}}(U_{s-}-Y_{s-})dA^{\tau *}_{s}+\sum_{\tau\leq s<\lambda_{\tau}}(U_{s-}-Y_{s-})\Delta^{+}A^{\tau}_{s}=0,\,\,\ \text{a.s.}
\end{array}
\right.
\end{multline}
Indeed, if $t\in[\delta_{\tau},\lambda_{\tau}]$, we have $K^{\tau}_{\lambda_{\tau}}-K^{\tau}_{t}=0$  and ~\eqref{eqtheorem2} is satisfied from ~\eqref{eqproposition5}. And if $t\in[\tau,\delta_{\tau}]$, then from ~\eqref{eqdefinition4} we have, $Y_{t}=Y_{\delta_{\tau}}+\int_{t}^{\delta_{\tau}}f(s,Y_{s})ds+K^{\tau}_{\delta_{\tau}}-K^{\tau}_{t}-\int_{t}^{\delta_{\tau}}dM^{\tau}_{s}.$
As $Y_{\delta_{\tau}}=Y_{\lambda_{\tau}}+\int_{\delta_{\tau}}^{\lambda_{\tau}}f(s,Y_{s})ds-(\tilde{A}^{\tau}_{\lambda_{\tau}}-\tilde{A}^{\tau}_{\delta_{\tau}})-\int_{\delta_{\tau}}^{\lambda_{\tau}}d\tilde{M}^{\tau}_{s},$
then ~\eqref{eqtheorem2} is also satisfied since $K^{\tau}_{\lambda_{\tau}}-K^{\tau}_{\delta_{\tau}}=0$.\\
Now for any $t\in[\tau,\lambda_{\tau}]$, $L_{t}\leq Y_{t}\leq U_{t}$  a.s. and
$\int_{\tau}^{\lambda_{\tau}}(Y_{s-}-L_{s-})dK^{\tau *}_{s}+\sum_{\tau\leq s<\lambda_{\tau}}(Y_{s-}-L_{s-})\Delta^{+}K^{\tau}_{s}=\int_{\tau}^{\delta_{\tau}}(Y_{s-}-L_{s-})d\bar{K}^{\tau *}_{s}+\sum_{\tau\leq s<\delta_{\tau}}(Y_{s-}-L_{s-})\Delta^{+}\bar{K}^{\tau}_{s}=0$ and
$\int_{\tau}^{\lambda_{\tau}}(U_{s-}-Y_{s-})dA^{\tau *}_{s}+\sum_{\tau\leq s<\lambda_{\tau}}(U_{s-}-Y_{s-})\Delta^{+}A^{\tau}_{s}=\int_{\delta_{\tau}}^{\lambda_{\tau}}(U_{s-}-Y_{s-})d\tilde{A}^{\tau *}_{s}+\sum_{\delta_{\tau}\leq s<\lambda_{\tau}}(U_{s-}-Y_{s-})\Delta^{+}\tilde{A}^{\tau}_{s}=0$. \\
Finally $Y_{T}=\xi$ and  the process $(Y_{t},M^{\tau}_{t},K^{\tau}_{t},A^{\tau}_{t})$  is a local solution for ~\eqref{eqdefinition4} on $[\tau,\lambda_{\tau}]$.
\end{proof}
\begin{remark}
\label{remark8}
\begin{itemize}
  \item If $M_{t}$  is a local martingale w.r.t $\mathcal{F}_{t}$ and if  $\tau$ and $\delta_{\tau}$ are  two $\mathcal{F}_{t}$-stopping times such that $\tau\leq\delta_{\tau}$, then $M_{t}1_{[\tau\leq t\leq\delta_{\tau}]}$ is a $\mathcal{F}_{t}$-martingale. Indeed:
\begin{multline}
M_{t}1_{[\tau,\delta_{\tau}]}(t)=\int_{0}^{t}1_{[\tau,\delta_{\tau}]}(s)dM_{s}+\int_{0}^{t}M_{s}d(1_{[\tau,\delta_{\tau}]}(s))
=\int_{0}^{t}1_{[\tau,\delta_{\tau}]}(s)dM_{s}-M_{\delta_{\tau}\wedge t}+M_{\tau\wedge t}
\end{multline}
 \item The construction of $Y$ does not depend on $\tau$ but the ones of $M$, $K$ and $A$ do.
\end{itemize}
\end{remark}
\section{Existence of a global solution for  reflected BSDE with
two completely separated barriers.}
We are now ready to give the main result of this paper. Let us assume that the barriers $L$ and $U$ and their left limits are completely separated, i.e., they satisfy the following assumption:
$$[\textbf{H}]:\,\,\,\ \mathbb{P}-a.s.,\,\,\,\ \forall t\in[0,T],\,\,\,\ L_{t}<U_{t} \,\,\,\ and \,\,\,\ L_{t-}<U_{t-}.$$
\begin{theorem}
\label{theorem3}
Under Assumption $[\textbf{H}]$, there exists a unique process $(Y_{t},M_{t},K_{t},A_{t})_{t\leq T}$ solution of the reflected BSDE associated with $(f,\xi,L,U)$. i.e., $Y$,$K$,$A\in\mathcal{S}^{2}$, $M\in\mathcal{M}_{loc}$ with $M_{0}=0$ and
\begin{multline}
\label{eqtheorem3}
\left\{
  \begin{array}{ll}
    Y_{t}=\xi+\int_{t}^{T}f(s,Y_{s})ds+K_{T}-K_{t}-(A_{T}-A_{t})-\int_{t}^{T}dM_{s}, \,\,\ \forall t\leq T \\
  L_{t}\leq Y_{t}\leq U_{t}, \,\,\ a.s. \,\,\ \forall t\leq T, \,\,\ \text{and}  \\
   \int_{0}^{T}(U_{s-}-Y_{s-})dA_{s}^{*}+\sum_{0\leq s<T}(U_{s-}-Y_{s-})\triangle^{+}A_{s}=\int_{0}^{T}(Y_{s-}-L_{s-})dK_{s}^{*}+\sum_{0\leq s<T}(Y_{s-}-L_{s-})\triangle^{+}K_{s}=0, \,\,\,\,\ \text{a.s.}
\end{array}
\right.
\end{multline}
\end{theorem}
\begin{proof}
Let $(Y_{t})_{t\leq T}$  be the  process defined in Theorem~\ref{theorem2}, then $L \leq Y \leq U$ and $Y_{T}=\xi$. Now let $(\tau_{n})_{n\geq 0}$  a sequence of stopping times such that $\tau_{0}=0$ and $\tau_{n+1}=inf\{t\geq \tau_{n},Y_{t}=U_{t}\}\wedge T$ and $\tau_{n+2}=inf\{t\geq \tau_{n+1},Y_{t}=L_{t}\}\wedge T$. Henceforth, for any $n\geq0$ there exists a triple $(M^{n}_{t},K^{n}_{t},A^{n}_{t})_{t\leq T}$ of  processes  such that the process $(Y_{t},M^{n}_{t},K^{n}_{t},A^{n}_{t})_{t\leq T}$ is a local solution for the reflected BSDE associated with $(f,\xi,L,U)$ on the set $[\tau_{n},\tau_{n+1}]$ (by Theorem~\ref{theorem2}).\\
By the same argument in  \cite[3.7. Theorem.]{HamadeneHassani2005} ( see also \cite[Theorem 4.1.]{HamadeneHassaniOuknine2010}, \cite[Theorem 4.1.]{Hassairi2016} or \cite[Theorem 4.1]{HamadeneHang2009}) we show that  $\mathbb{P}([\tau_{n}<T,\forall n\geq 0])=0$, $\mathbb{P}$-a.s. since $\mathbb{P}$-a.s., $\forall t$, $L_{t-}<U_{t-}$. Which means that for $\omega\in\Omega$ there exists $n_{0}(\omega)\geq 0$ such that $\tau_{n_{0}}(\omega)=T$. Next let us introduce the following processes $M$, $K$, $A$: $\mathbb{P}$-a.s., for any $t\leq T$, one sets
\begin{eqnarray*}
  K_{t} &=& K_{\tau_{n}}+(K^{n}_{t}-K^{n}_{\tau_{n}}) \,\,\,\,\  if \,\,\,\,\ t\in]\tau_{n},\tau_{n+1}]\,\,\,\ (K_{0}=0)\\
  A_{t} &=& A_{\tau_{n}}+(A^{n}_{t}-A^{n}_{\tau_{n}})\,\,\,\,\ if\,\,\,\,\ t\in]\tau_{n},\tau_{n+1}]\,\,\,\ (A_{0}=0) \\
  M_{t} &=& M_{t}1_{[0,\tau_{1}]}+\sum_{n\geq 1}M^{n}_{t}1_{]\tau_{n},\tau_{n+1}]}.\\
\end{eqnarray*}
Since the sequence $(\tau_{n})_{n\geq0}$ is $\mathbb{P}$-a.s. of stationary type and for any $n\geq0$,
$E([M]_{\tau_{n}})<\infty$ then $E([M]_{T})<\infty$,  $\mathbb{P}$-a.s..\\
Next let us show that $(Y,M,K,A)$ is the solution of the reflected BSDE$(\xi,f,L,U)$. For any $n\geq0$ we have:
$\mathbb{P}$-a.s. for all $t\in[\tau_{n},\tau_{n+1}]$,
\begin{equation}
\label{eqtheorem32}
Y_{t}=Y_{\tau_{n+1}}+\int_{t}^{\tau_{n+1}}f(s,Y_{s})ds+K_{\tau_{n+1}}-K_{t}-(A_{\tau_{n+1}}-A_{t})-\int_{t}^{\tau_{n+1}}dM_{s}.
\end{equation}
For any $n\geq0$ we have: $\mathbb{P}$-a.s.
$Y_{\tau_{n}}=Y_{\tau_{n+1}}+\int_{\tau_{n}}^{\tau_{n+1}}f(s,Y_{s})ds+K_{\tau_{n+1}}-K_{\tau_{n}}-(A_{\tau_{n+1}}-A_{\tau_{n}})-\int_{\tau_{n}}^{\tau_{n+1}}dM_{s}.$
Now for any $n\geq0$ and $m\geq n$ we have:
$Y_{\tau_{n}}=Y_{\tau_{m}}+\int_{\tau_{n}}^{\tau_{m}}f(s,Y_{s})ds+K_{\tau_{m}}-K_{\tau_{n}}-(A_{\tau_{m}}-A_{\tau_{n}})-\int_{\tau_{n}}^{\tau_{m}}dM_{s}.$
By the fact that $(\tau_{n})_{n\geq0}$ is of stationary type and  taking $m$ large enough we obtain: $\forall n \geq0$, $\mathbb{P}$-a.s.,

\begin{equation}
\label{eqtheorem33}
Y_{\tau_{n}}=\xi+\int_{\tau_{n}}^{T}f(s,Y_{s})ds+K_{T}-K_{\tau_{n}}-(A_{T}-A_{\tau_{n}})-\int_{\tau_{n}}^{T}dM_{s}.
\end{equation}
Now let $t\in[0,T]$ then there exists $n_{0}$ such that $t\in[\tau_{n_{0}},\tau_{n_{0}+1}]$. Using  ~\eqref{eqtheorem32} then ~\eqref{eqtheorem33} we obtain:
$Y_{t}=\xi+\int_{t}^{T}f(s,Y_{s})ds+K_{T}-K_{t}-(A_{T}-A_{t})-\int_{t}^{T}dM_{s},$
which means that  $(Y,M,K,A)$  verify equation $(\textbf{LU4})$  of Definition~\ref{definition3}. Finally the processes $K$, $A$ and $Y$ satisfy
$\int_{0}^{T}(Y_{s-}-L_{s-})dK_{s}^{*}+\sum_{0\leq s<T}(Y_{s-}-L_{s-})\triangle^{+}K_{s}=
\sum_{n\geq 0}(\int_{\tau_{n}}^{\tau_{n+1}}(Y_{s-}-L_{s-})dK^{n,*}_{s}+\sum_{\tau_{n}\leq s<\tau_{n+1}}(Y_{s-}-L_{s-})\triangle^{+}K^{n}_{s})=0$ by  definition of $K$ and ~\eqref{eqtheorem3}. In the same way we have  $\int_{0}^{T}(U_{s-}-Y_{s-})dA_{s}^{*}+\sum_{0\leq s<T}(U_{s-}-Y_{s-})\triangle^{+}A_{s}=0$. Then the process $(Y,M,K,A)$ is a
solution for the reflected BSDE$(\xi,f,L,U)$. Uniqueness is a direct consequence of the comparison theorem (Proposition~\ref{proposition1}).
\end{proof}
\section*{Acknowledgements}
I would like to thank the Editor for carefully handling the paper and the anonymous referees for their
valuable comments and suggestions for improving the quality of this work.

\end{document}